\documentclass[final]{siamltex}
\usepackage{graphicx}

\def \pio2{{\pi \over 2}}

\def \beq{\begin{equation}}
\def \eeq{\end{equation}}

\def \lf{\left(}
\def \rt{\right)}

\def \ep{\varepsilon}

\title{NONPERIODIC TRIGONOMETRIC POLYNOMIAL APPROXIMATION}
\author{Hillel Tal-Ezer\thanks{Academic College of Tel-Aviv Yaffo, Israel}}
\begin{document}
\maketitle
\begin{abstract}
The suitable basis functions for approximating periodic function are periodic, trigonometric functions. When the function is not periodic, a viable alternative is to consider polynomials as basis functions. In this paper we will point out the inadequacy of polynomial approximation and suggest to switch from powers of $x$ to powers of $\sin(px)$ where $p$ is a parameter which depends on the dimension of the approximating subspace. The new set does not suffer from the drawbacks of polynomial approximation and by using them one can approximate analytic functions with spectral accuracy. An important application of the new basis functions is related to numerical integration. A quadrature based on these functions results in higher accuracy compared to Legendre quadrature.
\end{abstract}
   

\pagestyle{myheadings}
\thispagestyle{plain}
\markboth{NONPERIODIC TRIGONOMETRIC POLYNOMIAL APPROXIMATION}{NONPERIODIC TRIGONOMETRIC POLYNOMIAL APPROXIMATION}

\section{Introduction }
Let us consider a finite sum approximation of a function  $f(x)$ 
\beq
f(x)\approx \sum_{k=0}^n a_k \psi_k(x),\qquad \ x\in [-L,L]. 
\eeq
If $f(x)$ and its derivatives are continuous and periodic then the suitable basis functions are
\beq
\psi_k(x)=e^{i\lf\frac{n}{2}-k\rt\frac{\pi x}{L}}
\eeq
and the error
goes to zero exponentially fast. If the function is not periodic, a spectral rate of convergence can be achieved  by using Chebyshev polynomials as basis functions. Approximating a function by finite sum of Chebyshev polynomials has drawbacks (see, for example~\cite{tref}). Despite the spectral rate of convergence, Chebyshev polynomials have peculiar characteristics and therefore do not provide an optimal set of basis functions. This peculiarity lies in the behavior of their derivatives. While Chebyshev polynomials have uniform behavior in the interval, their derivatives are non-uniform. In order to clarify this point let us assume that the interval is  $[-1,1]$. Chebyshev polynomials in this interval are defined as
\beq 
T_k(x)=\cos(k\arccos(x)).
 \eeq
 We have
 \beq
 \max_{x\in[-1,1]}|T_k(x)|=1
 \eeq
 and this maxima is achieved at $k-1$ points. Hence, $T_k(x)$ behaves uniformly. Let us consider now $\frac{dT_k}{dx}$.
 Upon defining
\beq
 x=\cos(\theta)
 \eeq
 we get
\beq
\frac{dT_k}{dx}=\frac{dT_k}{d\theta}\frac{d\theta}{dx}=\frac{k\sin (k\theta)}{\sin(\theta)}. 
\eeq
Hence, in the  vicinity of $x=0$ we have
\beq
\max |\frac{dT_k}{dx} |\approx k \qquad \qquad 
\eeq
while at the end points we have
\beq
|\frac{dT_k}{dx}(\pm 1)|=k^2 .
\eeq
Thus, there is no uniformity in the extrema values of the first derivative.
This discrepancy increases
for each additional derivative by a power of $ \ k,\ $  namely
\beq
\frac{|T_k^{(j)}(x_a)|}{|T_k^{(j)}(x_b)|} \approx k^j \label{3}
\eeq
where $x_a$ and $x_b$ are extrema points in the vicinity of $\ \pm 1\ $  
and $\ 0\ $
respectively. Approximating an analytic function by Chebychev expansion is
 highly
efficient when the approximated function exhibits a behavior similar
to ~(\ref{3}). As in the Fourier case, the efficiency deteriorates as the
 behavior of
the function approximated moves away from this pattern. A reasonable target should be to look for an approximation space where the basis functions and all
their derivatives exhibit uniform, or `almost' uniform,  behavior.

The peculiarity of Chebyshev approximation can be demonstrated also by considering interpolation. It is well known
that polynomial interpolation in equally distributed points
is not the right approach. The error function has large gradients at the
boundaries and in some cases, divergence can occur(Runge
 phenomenon)~\cite{davs}.
Insisting on polynomial interpolation, one should resort to non-uniform
 distribution
of points. A set of interpolating points which results in exponential rate
of convergence (for analytic functions) is
\beq x_i=\cos \left( \frac{i \pi}{n}  \right)\qquad i=0,\cdots , n . \label{5}
\eeq
We have
\beq
|\Delta x_{max}|=|\cos\lf\frac{\pi}{2}-\frac{\pi}{n}\rt|\approx \frac{\pi}{n}\label{60}
\eeq
\beq
|\Delta x_{min}|=|1-\cos\lf\frac{\pi}{n}\rt|\approx \frac{1}{2} \lf \frac{\pi}{n}\rt ^2\label{61}
\eeq
The uneven distribution of points does not make sense in the general case.
When the function
approximated does not have large gradients at the boundaries there is
no justification for concentrating points there. This inefficient strategy is
more pronounced when the large gradients of the function are away
from the boundaries. A way of quantifying the inefficiency of
Chebychev interpolation is in the fact that we need $\ \pi\ $
points per wavelength for resolution, compared to 2 points
per wavelength which is ``a law of nature" as expressed by
Nyquist criteria. It seems reasonable to look for a set of
basis functions
where the relevant interpolating points are evenly, or almost evenly,
distributed. 

In this paper we present a new set of  basis functions. The set is composed of nonperiodic trigonometric functions. It was first introduced in~\cite{kota}. The motivation in that case was to overcome the severe stability restriction which results from using Chebyshev polynomials for space discretization while solving time dependent pde's. Many researchers reported significant increases in efficiency while using the approximating method described in~\cite{kota} (e.g.~\cite{carc},~\cite{don},~\cite{javi},~\cite{jodi},~\cite{solo},~\cite{wais}). In the present paper we would like to deal with these basis functions just from approximation  viewpoint and to emphasize the advantages of using these functions as compared to polynomials.  Approximating a general function by a linear combination of these basis functions overcomes the drawbacks mentioned above and, for analytic functions, the approximation is spectrally accurate.
The paper is organized as follows. In Section 2 we describe the
proposed approximation subspace. The basis functions depend on a parameter $p,\ 0 < p < \frac{\pi}{2},\ $ and this parameter depends on $n$ ( the size of the subspace). 
In Section 3 we carry an analysis related to accuracy properties and in Section 4, analysis related to resolution properties. In Section 5 we use the new basis function for numerical integration.
 The paper is concluded in Section 6 in which we present numerical results.
\section{nonperiodic Trigonometric Polynomial Subspace}
\noindent
Without loss of generality we will consider functions in $[-1,1]$. Let $P_k(y)$ be a set of polynomials, orthogonal with respect to the inner product
\beq
<f,g>=\int \limits_{-1}^1 f(y)g(y)w(y)dy.
\eeq
By change of variables
\beq
y=\frac{\sin(px)}{\sin(p)}\label{778},
\eeq
where $p$ is a parameter in $\lf 0,\frac{\pi}{2}\rt$ ,
we get that
\beq
\psi_k(x)=P_k\lf \frac{\sin(px)}{\sin(p)}\rt\label{666}
\eeq
are orthogonal with respect to the inner product
\beq
<f,g>=\int \limits_{-1}^1 f(x)g(x)\bar w(x)dx
\eeq
where
\beq
\bar w(x)=w\lf \frac{\sin(px)}{\sin(p)}\rt \cos(px).
\eeq
Hence, the proposed new approximating subspace is 
\beq
S_n=\textrm{span}\{\psi_0(x),\dots,\psi_n(x)\}.\label{987}
\eeq

It is easily verified, by trigonometric identities, that $S_n$  can be written also as a subspace spanned by trigonometric polynomials. More precisely
\beq
S_n=\textrm{span}\{Q_0,\dots,Q_{n}\}
\eeq
where 
\beq
Q_k(x)=\cos(kpx) \qquad k=0,2,4,\dots ,\qquad  0 < p < \frac{\pi}{2} \label{111}
\eeq
 and
\beq
Q_k(x)=\sin(kpx) \qquad k=1,3,5,\dots , \qquad 0 < p < \frac{\pi}{2}. \label{222}
\eeq
A popular set of orthogonal polynomials is the set of Jacobbi polynomials $P_k^{(\alpha,\beta)}(y)$. These polynomials are orthogonal under the inner product
\beq
<f,g>=\int \limits_{-1}^1 f(y)g(y)w(y)dy \label{567}
\eeq
where
\beq
w(y)=\lf 1-y \rt ^\alpha\lf 1 + y\rt^\beta,\qquad -1< \alpha\  \textrm {and} \ -1 < \beta.
\eeq
Hence
\beq
\bar w(x)=\lf 1-\frac{\sin(px)}{\sin(p)}\rt ^\alpha\lf 1 + \frac{\sin(px)}{\sin(p)}\rt^\beta\cos(px).
\eeq
 Since 
\beq
\lim_{p\rightarrow 0}\frac{\sin(px)}{\sin(p)}=x,
\eeq
Jacobbi polynomials can be considered as $\psi_k$ functions in the extremal case when $p=0$. As will be shown in the next section, in general, the parameter $p$ should be close to the other extremal point, namely $\frac{\pi}{2}$.

Two important members of the Jacobbi polynomials family are Chebyshev and Legendre polynomials. Let $T_k(x)$ be Chebyshev polynomial then
\beq
\psi_k(x)=T_k\lf \frac{\sin(px)}{\sin(p)}\rt
\eeq
and the weight function is 
\beq
\bar w(x)=\frac {\cos(px)}{\sqrt{1-\frac{\sin^2(px)}{\sin^2(p)}}}.
\eeq
If $P_k(x)$ is Legendre polynomial then 
\beq
\psi_k(x)=P_k\lf \frac{\sin(px)}{\sin(p)}\rt
\eeq
and
\beq
\bar w(x)=\cos(px).
\eeq

\section{Approximating analytic functions by projection on $S_n$}
Let $f(x)$ be a  function continuous in $[-1,1]$.
Orthogonal projection of $f$ on $S_n$ results in
\beq
f_n(x)=\sum_{k=0}^n a_k \psi_k(x) \label{24}
\eeq
where
\beq
a_k=\frac{<f,\psi_k>}{<\psi_k,\psi_k>}.
\eeq
Since Chebyshev polynomials is the mostly used set of orthogonal polynomials, we will consider
\beq
\psi_k(x)=T_k\lf \frac{\sin(px)}{\sin(p)}\rt
\eeq
in the rest of this section. Almost all the theoretical results described here are relevant to any set of orthogonal polynomials.

Since ~(\ref{567}) then
\beq
 <\psi_k,\psi_k>=\frac{\sin(p)}{p}\int\limits_{-1}^1 \frac{ T_k^2(y)}{\sqrt{1-y^2}}dy=\beta_k\frac{\pi\sin(p)}{2p}
\eeq
where $\beta_0=2$ and $\beta_k=1, \ 1\le k$. Therefore
\beq
a_k=\frac{2}{\beta_k \pi}\int\limits_{-1}^1 \frac{\tilde f(y) T_k(y)}{\sqrt{1-y^2}}dy
\label{30}
\eeq
where
\beq
\tilde f(y)=f\lf g(y;p)\rt \label{31}
\eeq
and
\beq
g(y;p)=\frac{1}{p}\sin^{-1}\lf y\sin(p)\rt \label{90}
\eeq
is the inverse of~(\ref{778}).
Hence, approximating $f(x)$ 
by the new set of basis functions is equivalent to approximating
$\ f\lf g(y;p)\rt\ $ by Chebychev polynomials. As a result,  it is
sufficient to deal with the latter. Observe that now, due to the singularity of $\ g(y;p) $ at $\ y=\pm 1/\sin(p) \ $, the function approximated has large gradient at the boundaries and it is justified to use Chebyshev polynomials.
The relevant theory which discusses
polynomial approximation of functions with singularities
outside the domain of definition here follows~\cite{wals}.

Let  $K$  be a bounded continuum in  $C$  such that $\enskip K^c$, the
 complement
$K \enskip $ , is simply connected in the extended plane and contains the point
at infinity.
For such  $\ K\ $ there exist a conformal mapping
$\ \Psi (w)\ $  which maps the complement of
the unit disc onto $\ K^c\ $~\cite{wals}.
Let $\theta(y)$ be the inverse of $\Psi(w)$ and
\beq B_t=\{y:\ |\theta(y)|=t\} \qquad (\ t>1)  \label{38} \eeq
denote the level curves in $\ K^c\ $ then we have the following theorem:

\noindent
{\bf Theorem 2.1:}{\it Suppose $t>1$ is the largest number such that $F(y)$
is analytic inside $B_t$. The interpolating polynomials $P_n(y)$ with
interpolating points $\ y_i^n\ $ that are uniformly distributed on $K$
then satisfy
\beq
\lim_{n\rightarrow\infty}\max_{y\in K}|F(y)-P_n(y)|^{\frac{1}{n}}
=\frac{1}{t} . \label{390} \eeq}

\medskip
Since approximating an analytic function by Chebyshev polynomials is equivalent to interpolating the function at uniformly distributed points 
(e.g. Chebyshev points)
, the asymptotic rate of convergence can be computed by making use
of this theorem.
For $\ K=[-1,1],\ $  the relevant
conformal mapping
is~\cite{mark}
\beq \theta(y)=y\pm \sqrt{y^2-1}\ . \label{311} \eeq

 $\ \tilde f(y)\ $ is singular at $\ y=\pm 1/\sin(p) \ $ hence , the largest $t$ is

\beq
 t=\frac{1}{\sin(p)} + \sqrt{\frac{1}{\sin^2(p)}-1}=\frac{1+\cos(p)}{\sin(p)}=\cot\lf\frac{p}{ 2}\rt \label{313}
  \eeq
and the asymptotic rate of convergence is
\beq
\frac {1 }{ t}=\tan\lf\frac{p}{ 2}\rt. \label{312}
  \eeq
Hence, approximating by the new set of basis functions, the asymptotic accuracy is  $\ c\ep\ $ where
\beq
 \ep=\lf \tan\lf\frac{p}{ 2}\rt\rt^n  \label{700}
  \eeq
and  $\ c\ $ is  constant which depends on $\ \tilde f(y)\ $ but does not depend on
$\ n\ $ or $\ y.\ $

By choosing 
\beq
p=2\arctan\lf \ep^{\frac{1}{ n}}\rt, \label{991}
\eeq
where $\ep$ is the machine accuracy, we eliminate the error which results from the singular points $\ y=\pm 1/\sin(p) \ $ and get spectral accuracy. Detailed mathematical proof of spectral accuracy is given in~\cite{cost1}.

The choice of $p$ described above is independent of the function we are approximating. Obviously, one can do better by choosing the appropriate $p$ for each function. This can be achieved numerically by making use of a minimization algorithm which finds the parameter $p$ that minimizes the norm of the error vector defined as
\beq
E=\sum_{j=1}^m|\sum_{k=0}^n a_k\lf p\rt\psi_k(z_j)-f\lf z_j\rt | \label{990}
\eeq
where $z_j,\ 1\le j \le m\ ,$  are check points randomly distributed in the interval $[a,b]$.

\vspace{.1in}

Approximating $\ f(x) \ $ via interpolation, $\ a_k\ $ have to
satisfy the following $\ n+1\ $ equations
\beq
\sum_{k=0}^n a_k \psi_k(x_i)=f(x_i)\qquad 0\leq i\leq n  \label{91}
\eeq
where $\ \{x_i\}_{i=0}^n\ $ is an appropriate set of
interpolating points. Due to the equivalence mentioned above, a feasible
set of interpolating points is
\beq
x_i=g(y_i;p), \qquad
 0\leq i \leq n \label{95}
\eeq
where
\beq
 y_i=\cos\lf \frac{i\pi}{ n}\rt. \label{40}
\eeq

We would like  to show now that, while in the Chebyshev case we have
\beq
\lim_{n\rightarrow\infty} \frac{\Delta x_{min}}{\Delta x_{max}}=0,
\eeq
in the Nptp (nonperiodic trigonometric polynomial) case we have
\beq
\lim_{n\rightarrow\infty} \frac{\Delta x_{min}}{\Delta x_{max}}=c,\qquad c\ne 0.
\eeq

Let us address the general case 
\beq
\lim_{n\rightarrow\infty} \frac{\Delta x_i}{\Delta x_{max}}. \label{71}
\eeq
We have
\beq
\frac{\Delta x_i}{\Delta x_{max}}=\frac{\sin^{-1}\lf \sin(p)y_{i+1}\rt -\sin^{-1}\lf \sin(p)y_{i}\rt}{0-\sin^{-1}\lf \sin(p)y_{\frac{n}{2}-1}\rt}.
\eeq

Upon defining
\beq
\theta = \frac {\pi}{n}
\eeq
we have 
\beq
\frac{\Delta x_i}{\Delta x_{max}}=\frac{\sin^{-1}\lf \sin(p)\cos{\lf i\theta\rt}\rt -\sin^{-1}\lf \sin(p)\cos{\lf (i+1)\theta\rt}\rt}{\sin^{-1}\lf \sin(p)\sin(\theta)\rt}.
\eeq
Computing ~(\ref{71}) via l'Hopital's rule, we have to compute the derivatives of the numerator and denominator. 
(in what follows we will compute $\lim_{\theta \rightarrow 0}{}$ instead of  $\lim_{n \rightarrow \infty}{}$).

We have
\beq
\frac{d}{d\theta}\sin^{-1}\lf \sin(p)\cos{(i\theta)}\rt= \frac{-i\sin(p)\sin(i\theta)+\cos(i\theta)\cos(p)\frac{dp}{d \theta}}{\sqrt{1-\sin^2(p)\cos^2(i\theta)}}.
\eeq

One can write the r.h.s of the equation above as $h_1(\theta)h_2(\theta)$ where
\beq
h_1(\theta)=\frac{\sin(i\theta)}{\sqrt{1-\sin(p)\cos(i\theta)}\sqrt{1+\sin(p)\cos(i\theta)}}
\eeq
and
\beq
h_2(\theta)=-i\sin(p)+\frac{\cos(p)}{\sin(i\theta)}\cos(i\theta)\frac{dp}{d\theta}.
\eeq
Since (\ref{991}) we have
\beq
\frac{dp}{d\theta}=2\frac{1}{1+\varepsilon^{\frac{2\theta}{\pi}}}\varepsilon^{\frac{\theta}{\pi}}\ln\lf{\varepsilon}\rt\frac{1}{\pi}
\eeq
and therefore
\beq
\lim_{\theta \rightarrow 0}\frac{dp}{d\theta}=\mu
\eeq
where
\beq
\mu=\frac{\ln(\varepsilon)}{\pi}.
\eeq
Using l'Hopital's rule we get
\beq
\lim_{\theta \rightarrow 0}h_1\lf \theta\rt = \frac {i}{\sqrt{i^2+\mu^2}}
\eeq
and
\beq
\lim_{\theta \rightarrow 0}h_2\lf \theta\rt = -\frac {i^2+\mu^2}{i}.
\eeq
Hence
\beq
\lim_{\theta \rightarrow 0}\frac{d}{d\theta}\sin^{-1}\lf \sin(p)y_{i}\rt=-\sqrt{i^2+\mu^2}.
\eeq
As to the denominator, using  l'Hopital's rule again results in
\beq
\lim_{\theta \rightarrow 0}\frac{d}{d\theta}\lf \sin^{-1}\lf \sin(p)\sin(\theta)\rt\rt=1.
\eeq
Based on  the results above we finally get
\beq
\lim_{\theta \rightarrow 0}\frac{\Delta x_i}{\Delta x_{max}}=\sqrt{(i+1)^2+\mu^2}-\sqrt{i^2+\mu^2},\qquad  i\ge 0.\label{999}
\eeq
Since $\Delta x_0=\Delta x_{min}$ and using $\varepsilon=10^{-8}$, for example, we have $|\mu|=5.8635$ and therefore
\beq
\lim_{\theta \rightarrow 0}\frac{\Delta x_{min}}{\Delta x_{max}}=0.0847.
\eeq
Observing  ~(\ref{999}) we can conclude that the interpolating points are "almost" equally distributed as $n\rightarrow\infty$. For example,  
the number of points which satisfy
\beq
\lim_{\theta \rightarrow 0}\frac{\Delta x_i}{\Delta x_{max}} < 0.9 
\eeq
is only $22$. For $n$ large, this number is negligible.

Let us look now at the behavior of the derivatives of the basis functions compared to the Chebyshev case ~(\ref{3}).
Upon using~(\ref{666})
and defining
\beq
\cos\lf\theta\rt=\frac{\sin\lf px \rt }{ \sin\lf p \rt } \label{1000}
\eeq
we have
\beq
\psi_n(x)=\cos{n\theta}.
\eeq
Hence
\beq
 \frac{d\psi_n(x)}{ dx}=\frac{d\psi_n}{ d\theta}\frac{d\theta}{ dx} = n\frac{p}{\sin\lf p \rt }\frac{\sin(n\theta)}{\sin(\theta)}\cos(px). \label{09}
\eeq
The maxima of the derivative is achieved at $x=1 ( \theta=0)$ and the minima(for $n$ odd) at $x=0( \theta=\frac{\pi}{2})$ . Since ~(\ref{991}) then
\beq
\cos(p)=\frac{1-\varepsilon^{\frac{1}{n}}}{1+\varepsilon^{\frac{1}{n}}}
\eeq
and therefore, using l'Hopital's rule,
\beq
\lim_{n \rightarrow \infty}\frac{\max|\frac{d\psi_n(x)}{ dx}|}{\min|\frac{d\psi_n(x)}{ dx}|}=\frac{|\ln{(\varepsilon)}|}{2}.
\eeq
Hence, the nonuniformity of the first derivative almost diminishes. In a similar way one can show almost uniformity for higher derivatives.

\section{On Resolution}
 Let 
\beq
f(x)=\sin (r\pi x )\qquad -1\le x \le 1. \label{777}
\eeq
(similar analysis can be carried out for  $f(x)=\cos (r\pi x )$).

By change of variables $ y=\frac{\sin(px)}{ \sin(p)}$ we get
\beq
f(x)=\tilde f (y)
\eeq
where
\beq
\tilde f(y) = \sin\lf \frac{r\pi}{ p}\sin^{-1}\lf y\sin\lf p\rt \rt \rt .
\eeq
Hence, resolving $f(x)$ by projection on subspace spanned by $\psi_k(x)$, where
\beq
\psi_k(x)=T_k\lf \frac{\sin\lf px\rt}{\sin\lf p \rt }\rt,
\eeq
is equivalent to resolving $\tilde f (y)$ by Chebyshev polynomials. Let $r$ be chosen such that
\beq
 m=\frac{r\pi}{ p}\label{9090}
 \eeq
 is an odd number. We will show now that 
\beq
\tilde f(y)=(-1)^mT_m\lf \alpha y\rt
\eeq
where
\beq
\alpha=\sin\lf p \rt .
\eeq
\noindent
\underline{\bf Lemma 1:}

Let $T_k$ be Chebyshev polynomial of degree $k$ then, for $k$ even, we have
\beq
\cos(kpx)=(-1)^{\frac{k}{2}}T_k(\sin(px))
\eeq
and for $k$ odd
\beq
\sin(kpx)=(-1)^{\frac{k-1}{2}}T_k(\sin(px)) \label{500}
\eeq
\noindent
\underline{\bf Proof:}

The proof is by induction on $k$.  It is easily verified for $k=0,1$ .

Assume first that $k$ is even. By the recurrence relation of Chebyshev polynomials we get
\beq
T_{k+1}\lf\sin\lf px \rt \rt =2\sin\lf px \rt T_k \lf \sin \lf px \rt \rt -T_{k-1}\lf \sin\lf px \rt \rt
\eeq
and by induction 
\beq
T_{k+1}(\sin(px))=(-1)^{\frac{k}{2}}\lf 2\sin(px)\cos(kpx)+\sin((k-1)px)\rt.
\eeq
Since
\beq
2\sin(px)\cos(kpx)= \sin\lf \lf k+1 \rt px \rt -\sin \lf \lf k-1 \rt px \rt
\eeq
we get
\beq
T_{k+1}(\sin(px))=(-1)^{\frac{k}{2}}\sin((k+1)px)
\eeq
and the proof of the even case is concluded. 

For $k$ odd, using the recurrence relation and induction we have
\beq
T_{k+1}(\sin(px))=(-1)^{\frac{k-1}{2}}\lf 2\sin(px)\sin(kpx)-\cos((k-1)px)\rt.
\eeq
Since 
\beq
 2\sin(px)\sin(kpx)=\cos((k-1)px)-\cos((k+1)px)
 \eeq
then
\beq
T_{k+1}(\sin(px))=(-1)^{\frac{k+1}{2}}\cos((k+1)px)
\eeq
and the proof of the odd case is concluded.

Since $T_m\lf \alpha y\rt$ is polynomial of degree $m$ in $y$, it can be written as
\beq
T_m\lf \alpha y\rt=\sum_{k=0}^m  a_k^m T_k\lf y \rt 
\eeq
while
\beq
a_k^m=\frac{2}{\pi c_k}\int_{-1}^1 \frac{T_m\lf \alpha y\rt T_k\lf y \rt }{ \sqrt{ 1-y^2} } dy \qquad c_0=2, c_k=1 \ \textrm{for} \ k\ge 1.
\eeq
Chebyshev polynomials satisfy the recurrence relation

\beq
T_{n+1}(x)=2xT_n(x)-T_{n-1}(x). \label{600}
\eeq
Hence
\beq
a_k^m=\frac{2}{ \pi c_k}\int_{-1}^1 \frac{\left[ 2\alpha yT_{m-1}\lf\alpha y\rt-T_{m-2}\lf \alpha y\rt \right]T_k(y)}{ \sqrt{1-y^2}}dy\qquad k\ge0 ,\ m\ge 2.\label{300}
\eeq
Since ~(\ref{600}) we have 
\beq
2yT_{k}(y)=T_{k+1}(y)+T_{k-1}(y).
\eeq
Substituting this relation in ~(\ref{300}) we finally get that the coefficients satisfy
\beq
a_0^0=1,\ a_0^1 = 0, \ a_1^1=\alpha,
\eeq
\beq
a_0^m=\alpha a_1^{m-1}-a_0^{m-2} ,\qquad \ m\ge 2,
\eeq
\beq
a_k^m =  \alpha \lf c_{k-1}a_{k-1}^{m-1}+a_{k+1}^{m-1}\rt -a_k^{m-2} , \qquad 1\le k\le m , \qquad \ m\ge 2.
\eeq
Carrying out numerical experiments we have observed that , while $k < \alpha m$, $a_k^m$ oscillates, and once $k \ge \alpha m$, $|a_k^m|$ monotonically decreases. 
Based on this numerical results we had conjectured, in~\cite{kota}, that the function $T_m\lf \alpha y \rt $ is resolved by $k$ terms where $k =\lceil \alpha m \rceil$ ( it was proven later in~\cite{cost}).

Hence, by using $ \{\psi_k(x)\}_{k=0}^n$, the maximal $k$ is $n$ and therefore, since ~(\ref{9090}), we get
\beq
r  < r_{\max}
\eeq
where
\beq
r_{\max}=\frac{np}{ \pi \sin(p)}.
\eeq 
Hence
\beq
\lim_{p\rightarrow \frac{\pi}{ 2}}r_{\max}=\frac{n}{ 2}
\eeq
which is exactly Nyquist criteria.

\section{Numerical Integration}
Approximating 
\beq
I=\int_{-1}^1 f(x)dx \label{444}
\eeq
is an essential subject in numerical analysis. An highly accurate approach is Gaussian quadrature based on Legendre polynomials. 
In ~\cite{hale}, the authors describe nonpolynomial algorithms which are aimed at overcoming the "waste" of factor $\frac{\pi}{2}$ typical to polynomial algorithms.
A quadrature based on $\psi_k(x)$ is a member of this family of nonpolynomial algorithms. 

Using quadrature based on 
\beq
\psi_k(x)=P_k\lf \frac{\sin(px)}{\sin(p)}\rt,
\eeq
where $P_k$ is Legendre polynomial, is equivalent to doing first change of variables
\beq
y=\frac{\sin(px)}{\sin(p)}
\eeq
which transforms  ~(\ref{444}) to
\beq
I= \int_{-1}^1 \tilde f(y)dy \label{445}
\eeq
where
\beq
\tilde f(y)=\frac{\sin(p)}{p}\frac{f\lf \frac{1}{p}\sin^{-1}\lf y\sin(p)\rt\rt}{\sqrt{1-\lf\ y\sin(p)\rt^2}}
\eeq
and approximating ~(\ref{445}) by standard Legendre Gaussian quadrature. Hence the quadrature can be written as
\beq
\tilde I =\sum_{i=1}^m f(x_i)w_i
\eeq
where 
\beq
x_i=\frac{1}{p}\sin^{-1}\lf y_i\sin(p)\rt\label{800}
\eeq
while $\{y_i\}_{i=1}^m$ are the zeros of Legendre polynomial of degree $m$. $\{w_i\}_{i=1}^m$ are the weights defined as 
\beq
w_i=\frac{\sin(p)}{p}\frac{\bar w_i}{\cos\lf p x_i\rt}\label{801}
\eeq
while $\bar w_i$ are Legendre weights. Similarly one can write a quadrature which is related to Chebyshev polynomials or any other set of orthogonal polynomials.

\section{Numerical Results}

\noindent
In the first set of examples we present results related to approximating functions. In this set we are comparing two algorithms :

1. Chebyshev

2. Nptp (Non periodic trigonometric polynomial).

The error presented in the tables  is defined as 
\beq
Er = \sqrt { \sum _{j=1}^{100} \lf f(y_j) - \tilde f (y_j) \rt ^2 }
\eeq
where $f$ is the exact  function we are approximating, $\tilde f $ is the approximating function which results from using either Chebyshev polynomials or Nptp
functions and 
\beq
y_j=a+(j-1)\frac{b-a}{ 99},\qquad 1\le j \le 100
\eeq
are check points, equally distributed in the interval $[a,b]$.

For Nptp, the tables contain two sets of results. One set, Nptp1, is for the case where the parameter $p$ is computed according to ~(\ref{991})   
\beq
p=2\arctan(\ep^{\frac{1}{ n}}),\qquad \ep=10^{-15}
\eeq
and the second, Nptp2, is for the case where the parameter $p$ is computed adaptively by minimizing ~(\ref{990}).  The number in the brackets contains the parameter $p$ in each case.  

\vspace{.2in}

\noindent
\underline{\bf Example 1}

In this example we approximated 

\beq
f(x)=\frac{1}{ 2+\cos(40x)},\qquad -1\le x\le 1.
\eeq
This function behaves uniformly. The results are
\vspace{.3in}

\begin{tabular}{|c|c|c|c|}
\hline
n & Er(Chebyshev) & Er(Nptp1) & Er(Nptp2)\\
\hline
100 &4.7562e-2&1.5344e-2 (1.232)&1.4528e-2 (1.433) \\  
200 &2.2647e-3&7.6117e-5 (1.399)& 5.1861e-5 (1.468)\\
400& 2.8352e-6&7.9950e-9 (1.485)&7.8335e-9 (1.390) \\ 
\hline 
\end{tabular} 

\vspace{.2in}
\noindent
\underline{\bf Example 2}
In this example we approximated the function 
\beq
f(x)=x^5\cos\lf 50x\rt,\qquad -1\le x\le 1.
\eeq
Due to the $x^5$ term, the function has large gradients close to the boundaries, nevertheless, the new algorithm outperforms Chebyshev approximation as can be seen by the results presented in the next table

\vspace{.3in}

\begin{tabular}{|c|c|c|c|}
\hline
n & Er(Chebyshev) & Er(Nptp1) & Er(Nptp2)\\
\hline
40 &6.0165e-2&3.5717e-2 (0.840)&3.7053e-2 (0.852) \\  
50 &3.7617e-3&5.4146-4 (0.967)& 2.0321e-6 (1.086)\\
60&1.6279e-4&4.5186e-11 (1.058)&3.3845e-11 (1.057) \\ 
\hline 
\end{tabular}

\vspace{.2in}

\noindent
\underline{\bf Example 3}

In this example we approximated the function 
\beq
f(x)=e^{-30x^2},\qquad -1\le x\le 1
\eeq
 and the results are
   \vspace{.3in}

\begin{tabular}{|c|c|c|c|}
\hline
n & Er(Chebyshev) & Er(Nptp1) & Er(Nptp2) \\
\hline
10 & 4.8234e-1&4.8138e-1 (0.0796)&1.3220e-1 (1.5708)\\  
20 & 2.9417e-2&2.4545e-2 (0.3939)&2.0958e-4 (1.5708) \\
40 & 2.9475e-6&7.3752e-8 (0.8402)&4.5169e-14 (1.5708) \\ 
\hline 
\end{tabular}
\vspace{.2in}

In this case, due to the fact that there is a large gradient in the center of the interval and that the function is almost $0$ at the boundaries, the accuracy is significantly improved by choosing the optimal parameter (Nptp2) which, in this case, is very close to $\frac{\pi}{ 2}$. 
\vspace{.2in}

\noindent
\underline{\bf Example 4}

In this example we approximated the function 
\beq
f(x)=\frac{1}{ \sqrt{1.1-x^2}},\qquad -1\le x\le 1
\eeq
 and the results are
   \vspace{.3in}

\begin{tabular}{|c|c|c|c|}
\hline
n & Er(Chebyshev) & Er(Nptp1) & Er(Nptp2) \\
\hline
20 & 2.5252e-3 &2.8448e-3 (0.3939)   &2.5252e-3 (0) \\  
40 &3.7085e-6 &2.0681e-5 (0.8402)& 3.7085e-6 (0) \\
80 & 9.7418e-12 &3.3488e-8 (1.1783) & 9.7418e-12 (0)\\ 
\hline 
\end{tabular}
\vspace{.2in}

In this case , the behavior of the function at the boundaries justifies interpolating at Chebyshev points. As expected, the minimization process resulted with $p=0$ which means choosing Chebyshev points.

\vspace{.2in}
\noindent
\underline{\bf Example 5}

In this example we would like to demonstrate the resolution properties of Nptp compared to Chebyshev. For this purpose we approximated  the function
\beq
f(x)=\sin(100\pi x)+\cos(100\pi x),\qquad -1\le x\le 1.
\eeq
The results are presented in the next two tables.

\vspace{.2in}

\begin{tabular}{|c|c|c|}
\hline
n & p & Er(Nptp1) \\
\hline
220 &  1.4248 & 3.5265e-1 \\
240 &  1.4369 & 4.9448e-7 \\ 
260 &   1.4471 &3.6805e-9 \\ 
\hline 
\end{tabular}

\vspace{.2in}

\begin{tabular}{|c|c|}
\hline
n &  Er(Chebyshev)\\
\hline
320  &  2.4532e-1 \\
340 &    2.6849e-4 \\ 
360 &   1.6117e-8 \\ 
\hline 
\end{tabular}

\vspace{.2in}

In the next set of examples we use Nptp to approximate definite integrals
\beq
I(f)=\int_{-1}^1 f(x)dx.
\eeq

The tables below present the absolute value of the error while using 2 methods: Legendre and Nptp. 

In the first table, the function is
\beq
f(x)=\frac{100\cos\lf 100 x\rt}{2+\sin\lf 100 x \rt}.
\eeq
The parameter $p$ is computed by ~(\ref{990}) with $\varepsilon=1.e^{-5}$.

\vspace{.2in}

\begin{tabular}{|c|c|c|c|}
\hline
n & ErLegendre & ErNptp\\
\hline
200 &  6.2532e-2 & 1.0331e-3 \\
300 &  4.5825e-3 & 3.7822e-6  \\ 
500 &   1.2392e-5 & 1.8049e-9 \\ 
\hline 
\end{tabular}

\vspace{.2in}

In the next table the function is
\beq
f(x)=\cos(500x)
\eeq
 and $p$ is computed by ~(\ref{990}) with $\varepsilon=1.e^{-15}$.
 \vspace{.2in}

\begin{tabular}{|c|c|c|c|}
\hline
n & ErLegendre & ErNptp\\
\hline
180 & 1.9069e-1 & 1.8320e-2\\
190 & 7.3531e-2 & 1.6238e-11 \\
200 & 2.2017e-1 & 2.0517e-14  \\
250 &  3.1385e-1 &3.0422e-14  \\ 
270 &  3.0560e-6 &1.0923e-14  \\ 
290 & 7.3459e-15   &  8.1304e-15 \\ 

\hline 
\end{tabular}

 \vspace{.2in}

 Observe that, while in the Nptp case $200$ points were enough to recover the solution with machine accuracy, in the standard Legendre quadrature we needed almost $50\%$ more points in order to get machine accuracy.  
\vspace{.2in}

\noindent
\underline{\bf Conclusions:} We have presented in this paper a new set of basis functions which can be used for approximating general, smooth function defined on a real interval $[a,b]$. The new space is spanned by powers of trigonometric functions instead of powers of $x$ as in the regular polynomial case.
The trigonometric functions depend on a parameter $p$ which is a function of the dimension of the approximating subspace. When one fixes $p$ to be zero he gets polynomials. Hence, polynomials can be considered as a singular member of the family where $p$ is fixed and equal to zero. As described in the paper and verified by numerical experiments, the parameter $p$ should approaches the other extreme value, $\frac{\pi }{ 2}$, as the dimension increases. Besides exponential accuracy,  the approximating function can be computed efficiently using FFT.  Approximating a function by the new set is equivalent to approximating a transformed function by polynomials. Hence, the vast literature related to polynomials can be used for analyzing algorithms which make use of the proposed set of functions.

\end{document}